\magnification=\magstep1
\vbadness=10000
\hbadness=10000
\tolerance=10000

\proclaim Coxeter groups, Lorentzian lattices, and K3 surfaces.
\hfill 25 June, 27 July 1998

Richard E. Borcherds, 
\footnote{$^*$}{ Supported by a Royal Society
professorship and an NSF grant.}

D.P.M.M.S.,
16 Mill Lane,
Cambridge,
CB2 1SB,
England.

e-mail: reb@dpmms.cam.ac.uk

www home page www.dpmms.cam.ac.uk/\~{}reb

\bigskip

\proclaim Contents.

1.~Introduction.

2.~Notation and statement of main theorem.

3.~Proof of main theorem.

4.~The structure of $\Gamma_{\Omega}$.

5.~Examples.

\proclaim
1.~Introduction.

The main result of this paper describes the normalizer
$N_{W_\Pi}(W_J)$ of a finite parabolic subgroup $W_J$ of a (possibly
infinite) Coxeter group $W_\Pi$.  More generally we describe
$N_{W_\Pi.\Gamma_\Pi}(W_J)$ where $\Gamma_\Pi$ is a group of diagram
automorphisms of the Coxeter diagram $\Pi$ of $W_\Pi$. By taking $\Pi$
to be Conway's Coxeter diagram of the reflection group of $II_{1,25}$
we compute the automorphism groups of some Lorentzian lattices and K3
surfaces.

In the case when $W_\Pi$ is a finite Coxeter group (and $\Gamma_\Pi =
1$) the normalizer of $W_J$ has been described by Howlett [H]. His
result states that $N_{W_\Pi}(W_J)$ is a split extension $W_J.W'_J$,
where $W'_J=W_{\Omega}.\Gamma_{\Omega}$ is in turn a split extension
with $W_{\Omega}$ a Coxeter group and $\Gamma_{\Omega}$ a more
mysterious group acting on $\Omega$. Howlett showed by case by case
analysis that if $\Pi$ is connected (and $W_\Pi$ is finite) then
$\Gamma_{\Omega}$ is an elementary abelian  2-group
and is a subgroup of $Aut(J)$.
When $W_\Pi$ is infinite the normalizer $N_{W_\Pi}(W_J)$ has a similar
structure, except that the group $\Gamma_{\Omega}$ can be more
complicated. Although there is still a canonical map from
$\Gamma_\Omega$ to $Aut(J)$, the kernel can be non trivial, though it
has finite cohomological dimension. The kernel is trivial in the case
of finite $W_\Pi$ considered by Howlett because any finite group of
finite cohomological dimension must be trivial.  For example, the case
when $J=A_1$, $\Gamma_\Pi=1$ has been done by Brink [Br], who showed
that $\Gamma_\Omega$ is a free group (and therefore has cohomological
dimension at most 1). We will extend Brink's result to Coxeter
diagrams of arbitrary finite reflection groups. More precisely we
construct a category $Q_4$ using $\Pi$ and $J$ and prove that the
classifying space of this category is a classifying space of the group
$\Gamma_{\Omega}$. The main point about this category $Q_4$ is that it is
often finite and can often be written down explicitly,
in which case we can easily read off a presentation of $\Gamma_\Pi$.
For example, if $J=A_1$ we show that the
classifying space of this category $Q_4$ is 1-dimensional, so its
fundamental group is free and we recover Brink's result.
After writing this paper I discovered that Brink and Howlett
had previously announced a related description
of the normalizer of a parabolic
subgroup of a Coxeter group; see [B-H], and
see example 2.8 for the relation between their
result and theorem 2.7.

For later applications we need some generalizations as follows. First
of all, instead of calculating the normalizer of $W_J$ in a Coxeter group
$W_\Pi$, we calculate the normalizer in an extension
$W_\Pi.\Gamma_\Pi$, where $\Gamma_\Pi$ is a group of diagram
automorphisms. Secondly, we sometimes want to compute not the full
normalizer, but a subgroup with image contained in some subgroup
$\Gamma_J$ of $Aut(J)$. Thirdly, we sometimes want to vary the choice
of the Coxeter group $W_{\Omega}$, which we do by varying
a certain normal subgroup $R$ of $\Gamma_J$. For example, in calculating the
automorphism groups of $K3$ surfaces we take $W_{\Omega}$ to be
generated by reflections of norm $-2$ vectors rather than by all
reflections, so we take $R=1$.

Section 2 contains a statement of the main result (theorem 2.7)
describing a classifying category for the group $\Gamma_\Omega$,
and section 3 contains the proof of this result. Section 4 contains 
some more information about the structure of $\Gamma_\Omega$. 
In section 5 we give some applications of theorem 2.7, and
in particular
show how to describe the automorphism groups of some Lorentzian lattices
by embedding them in $II_{1,25}$ and using the description of
$Aut(II_{1,25})$ in [C]. The idea of studying Lorentzian lattices
by embedding them as orthogonal complements of root lattice in
$II_{1,25}$ comes from Conway and Sloane ([C-S]).

Work of  I. Piatetski-Shapiro and I. R. Shafarevich [P-S]
shows that there is a map from the automorphism  group of a K3 surface
to the group of automorphisms of its Picard lattice
modulo the group generated by reflections of norm $-2$ vectors
which has finite kernel and co-finite image,
so in practice if we want to describe the automorphisms of a K3 
surface the main step is to calculate the automorphism  group of its
Picard lattice. 
Kondo showed in [K] that the automorphism groups of some
K3 surfaces could be studied by embedding their Picard lattice
as the orthogonal complement of a root lattice in $II_{1,25}$.
We use Kondo's idea to describe the automorphism groups
of some K3 surfaces in terms of combinatorics of the Leech lattice. 
In particular we reprove some results of
Vinberg [V] on the ``most algebraic'' K3 surfaces and extend them to the
``next most algebraic'' K3 surface. Kondo showed in [K] that the
automorphism group  of the Kummer surface of a generic genus 2 Jacobian 
was generated by the classically known automorphisms together with some
new automorphisms found by Keum [Ke], and we show how to use Kondo's results
to describe the structure of this group. Kondo and Keum [K-K] have recently
proved similar results for some Kummer surfaces associated to the products
of two elliptic curves. 

Kondo recently found another mysterious connection between 
automorphism  groups of K3 surfaces and Niemeier lattices [K98], 
and used this to give a short proof of Mukai's classification [Mu] of
the finite groups that act on K3 surfaces. 

I would like to thank
I. Cherednik, I. Grojnowski, R. B.  Howlett, 
J. M. E. Hyland, S. Kondo, U. Ray, and G. Segal
for their help.

\proclaim
2.~Notation and statement of main theorem.

This section states the main result (theorem 2.7) describing normalizers of
parabolic subgroups of Coxeter groups.

We recall some basic definitions about Coxeter systems.
For more about them see [Hi] or [Bo].
A pair $(W,S)$ is called a {\bf Coxeter system} if
$W$ is a group with a subset $S$ such that $W$ has the presentation
$$\langle s:s\in S|(ss')^{m_{ss'}}=1 \hbox{ when $m_{ss'}<\infty$}\rangle$$ 
where $m_{ss'}\in
\{1,2,3,\ldots,\infty\}$ is the order of $ss'$, and $m_{ss'}=1$ if and
only if $s=s'$. A diagram automorphism of $S$ is an automorphism of
the set $S$ that extends to an automorphism of the group $W$, and
$Aut(S)$ means the group of diagram automorphisms of $S$.  We say that
$(W,S)$ is {\bf irreducible} if $S$ is not a union of two disjoint
commuting subsets. The number of elements of $S$ is called its {\bf
rank}. The Coxeter system is called {\bf spherical} if $W$ has finite
order. The irreducible spherical Coxeter diagrams are $A_{n}$ $(n\ge
1)$, $B_n=C_n$ $(n\ge 2)$, $D_{n}$ $(n\ge 4)$, $E_6$, $E_7$, $E_8$,
$F_4$, $G^{(n)}_2=I_2(n)$ $(n\ge 5)$, $H_3$, and $H_4$. It is also
sometimes useful to define the Coxeter diagrams $B_1=C_1=A_1$,
$D_3=A_3$, $D_2=A_1^2$, $E_5=D_5$, $E_4=A_4$, $E_3=A_2A_1$,
$G_2^{(4)}=B_2=C_2$, $G_2^{(3)}=A_2$, $G_2^{(2)}=A_1^2$.

If $(W_\Pi,\Pi)$ is a Coxeter system then we write
$V_\Pi$ for the (possibly infinite dimensional) real 
vector space with a basis of elements $e_s$
for $s\in \Pi$, and put a symmetric bilinear form on $V_\Pi$ by
defining
$$(e_s,e_{s'})= 2\cos(\pi/m_{ss'}).$$
Note that we normalize the roots $e_s$ so that they have norm $(e_s,e_s)=-2$
rather than 1; this is done to be consistent with the usual conventions
in algebraic geometry.

The Coxeter group $W_\Pi$ acts
on $V_\Pi$ with the element $s\in \Pi\subseteq W_\Pi$ acting as 
the reflection $v\mapsto v+(v,e_s)e_s$ in the
hyperplane $e_s^\perp$. Any subgroup $\Gamma_\Pi$ of $Aut(\Pi)$
acts on $V_\Pi$ by
permutations of the elements $e_s$, so we get an action of
$W.\Gamma_\Pi$ on $V_\Pi$, and hence on the dual space $V_\Pi^*$. 
We write $\Delta^+$ for the set of positive roots of $W$.
We define the fundamental
chamber $C_\Pi\subseteq V_\Pi^*$ of $W_\Pi$ by
$$C_\Pi=\{x\in V_\Pi^*|x(r)\ge 0 \hbox{ for all } r\in \Pi
\hbox{ and } x(r) > 0 \hbox{ for almost all } 
r\in \Delta^+ \}.$$
(Recall that ``almost all'' means ``all but a finite number of''.)
A theorem due
independently to Tits and Vinberg
states that no two distinct points of $C_\Pi$ are conjugate under
$W$, and the subgroup of $W$ fixing all points of some subset $A$ of
$W$ is generated by the reflections in the faces of $W$ containing
$A$. In particular $W$ acts simply transitively on the conjugates
of $C_\Pi$. The union $W_\Pi(C_\Pi)$ of all conjugates of $C_\Pi$
under $W_\Pi$ is given by
$$W_\Pi(C_\Pi)=\{x\in V_\Pi^*| x(r) > 0 \hbox{ for almost all 
 } r\in \Delta^+ \}.$$ In particular $W_\Pi(C_\Pi)$ is
convex and closed under multiplication by positive real numbers.
If $x\in W_\Pi(C_\Pi)$ then the set of roots vanishing on $x$
is a finite root system.

Note that $W_\Pi(C_\Pi)$ is usually slightly smaller than the Tits cone,
which is defined in the same way except that we omit the condition
that $x(e_s) > 0$ for all but a finite number of $ s$ in the definition
of the fundamental domain. The Tits cone can be thought of
as obtained from $W_\Pi(C_\Pi)$ by ``adding some boundary components''.
The reason for
using $W_\Pi(C_\Pi)$ rather than the Tits cone is that
the cone $W_\Pi(C_\Pi)$ has the property that the subgroup
of $W_\Pi$ fixing any vector of it is finite.

We fix a spherical subset $J$ of $\Pi$.
In particular we get a
spherical Coxeter system $(W_J,J)$.
Suppose $K$ is an isometry of Coxeter diagrams
from $J$ into $\Pi$.
We write $W_K$ for the finite reflection group
generated by $K(J)$.
There is a natural homomorphism
$p:N_{W_\Pi.\Gamma_\Pi}(W_K)\mapsto Aut(J)$. We let
$N_{W_\Pi.\Gamma_\Pi}(W_K;\Gamma_J)$ be the subgroup
of elements whose image is in a subgroup $\Gamma_J$ of $Aut(J)$.
We are interested in describing the group
$N_{W_\Pi.\Gamma_\Pi}(W_J;\Gamma_J)$.
Most of the time we take $\Gamma_J=Aut(J)$
in which case $N_{W_\Pi.\Gamma_\Pi}(W_J;\Gamma_J) =N_{W_\Pi.\Gamma_\Pi}(W_J)$,
but it is occasionally useful to use other values of $\Gamma_J$;
see example 5.7 and theorem 4.1.

We define $W'_K$ to be the subgroup of $N_{W_\Pi.\Gamma_\Pi}(W_K;\Gamma_J)$
mapping $K(J)$ to itself.

\proclaim Lemma 2.1. $N_{W_\Pi.\Gamma_\Pi}(W_K;\Gamma_J) = W_K.W'_K$.

Proof. This follows immediately from the fact that
$W_K$ acts simply transitively on the  Weyl chambers of $W_K$, and
$W'_K$ is the subgroup of $N_{W_\Pi.\Gamma_\Pi}(W_K;\Gamma_J)$
fixing a Weyl chamber of $W_K$.
This proves lemma 2.1.

The group $W'_K$ acts on the subspace $V_\Pi^{*W_K}$ of $V_\Pi^*$ of
all vectors fixed by $W_K$. We now construct a reflection group
$W_{\Omega_K}$ acting on $V_\Pi^{*W_K}$.  We choose a normal subgroup
$R$ of $\Gamma_J$.  (The subgroup $R$ is used to control the
reflection group $W_{\Omega}$ defined below.  Often we want
$W_{\Omega}$ to be as large as possible and we take $R=\Gamma_J$, but
sometimes we want to take a smaller $W_{\Omega}$; see examples 5.3,
5.4 and 5.5.)  We define the group $W_{\Omega_K}$ to be the subgroup
of $W_\Pi\cap W'_K$ generated by elements $w\in W_\Pi\cap W'_K$ such
that $w$ acts on $V_\Pi^{*W_K}$ as a reflection and acts on $J$ as an
element of $R$.  We define $\Omega_K$ to be the Coxeter diagram of
$W_{\Omega_K}$ and $C_{\Omega_K}$ to be its fundamental chamber. If
$K$ is the identity map from $J$ to $J\subseteq \Pi$ then we write
$W_\Omega$, $C_\Omega$, and $\Omega$ instead of $W_{\Omega_K}$,
$C_{\Omega_K}$, and $\Omega_K$. Note that $W_{\Omega_K}$ is obviously
contained in the inverse image of $R$ in $W_\Pi\cap W'_K$, but can be
much smaller; see for example the discussion of $D_4$ in example 5.7.

We define the group $\Gamma_{\Omega_K}$ to be the subgroup of $W'_K$
of elements $w$ with $w(C_{\Omega_K})=C_{\Omega_K}$.

\proclaim Lemma 2.2. The group $W'_K$ is
a semidirect product $W'_K=W_{\Omega_K}.\Gamma_{\Omega_K}$.

Proof.
We first show that the group $W_{\Omega_K}$ acts faithfully on $V_\Pi^{*W_K}$.
 More generally we will show that if $w\in W_\Pi\cap W'_K$ acts
trivially on $V_\Pi^{*W_K}$, then $w=1$. To see this we
observe that $w$ fixes the point $ x\in C_\Pi\cap
V_\Pi^{*W_K}$ such that $x(e_s)=0$ if $s\in K(J)$ and $x(e_s)=1$ if 
$s\notin K(J)$.
Therefore $w$ is in the subgroup $W_K$ of $W_\Pi$ generated by the
simple reflections of $W_\Pi$ fixing $x$. On the other hand
$w$ maps $K(J)$ into itself as $w\in W'_K$.
This implies that $w=1$ because $1$ is the only element of
$W_K$ mapping $K(J)$ into itself. This
proves that the group $W_{\Omega_K}$ acts faithfully on $V_\Pi^{*W_K}$.

Lemma 2.2 now follows from the fact that $W_{\Omega_K}$ acts simply
transitively on the conjugates of $C_{\Omega_K}$ under $W'_K$, and
$\Gamma_{\Omega_K}$ is the stabilizer of $C_{\Omega_K}$. This proves
lemma 2.2.

Warning: the group $\Gamma_{\Omega_K}$ need not act faithfully on
$V_\Pi^{*W_K}$ (though it does act faithfully on $V_\Pi^{*W_K}\times J$).

We define a {\bf classifying category} of a group $\Gamma$ to be a
category whose geometric realization is a classifying space for
$\Gamma$. (Recall from [Q] that the geometric realization of
a category is a space with a 0-cell for each object and an $n$-cell for
each sequence of $n$ composable morphisms if $n>0$.)
For example, the category with one object whose morphisms
are the elements of $\Gamma$ (with composition given by group
multiplication) is a classifying category for $\Gamma$.

We have more or less reduced the problem of describing
$N_{W_\Pi.\Gamma_\Pi}(W_J; \Gamma_J)$ to that of describing
$\Gamma_{\Omega_K}$. (The Coxeter diagram $\Omega$ of $W_\Omega$
can be described once we know $\Gamma_\Omega$.)
The main theorem of this paper describes
$\Gamma_{\Omega_K}$ by giving an explicit classifying category for
it. To define this category we need some more definitions.

Suppose that $S$ is the Coxeter diagram of a finite reflection group
$G$ of a finite dimensional vector space with no vectors fixed by $G$.
Fix a Weyl chamber $C$ of $G$, so that the walls of $C$ correspond to
the points of $S$. Then there is a unique element $\sigma_S$ of $G$
taking $C$ to $-C$ called the opposition involution. The
involution $-\sigma_S$ acts on the the Coxeter diagram $S$, and its
action on $S$ does not depend on the choice of $C$. This action can be
described as follows. The points of the Coxeter diagram correspond to
the simple roots of $C$. This set of roots is the same as the set of
simple roots of $\sigma_S(C)=-C$ multiplied by $-1$.
Hence $-\sigma_S$ acts on this set of
simple roots, in other words on the Coxeter diagram $S$. The involution
$-\sigma_S$ of $S$ can be described explicitly as follows. On diagrams
of type $A_1$, $B_n=C_n$, $D_{2n}$, $E_7$, $E_8$, $F_4$, $G^{(2n)}_2$,
$H_3$, $H_4$ for any $n\ge 1$ the involution $-\sigma_S$ is
the trivial automorphism of $S$, while for diagrams of types
$A_{n+1}$, $D_{2n+1}$, $E_6$, $G_2^{(2n+1)}$ for $n\ge 1$ the involution
$-\sigma_S$ is the unique nontrivial automorphism of the Coxeter
diagram $S$. Finally if the diagram $S$ is a union of connected
components then $-\sigma_S$ acts on each connected component as
described above.

Suppose that $J$ and $S$ are Coxeter diagrams.  Suppose that $K$ and
$K'$ are two isometries from $J$ into $S$.  We define $K$ and $K'$ to
be {\bf adjacent} if there is a point $s$ of $S$ not in $K(J)$ such that
$K(J)\cup s$ is spherical and $\sigma_{K(J)\cup s}\sigma_{K(J)}$ takes
$K$ to $K'$. If $K$ is adjacent to $K'$ then $K'$ is adjacent to
$K$. We define two isometries $K$, $K'$ of $J$ into $S$ to be {\bf
associate} if there is a sequence of isometries $K=K_1$, $K_2,\ldots,
K_n=K'$ from $J$ into $S$ such that $K_i$ and $K_{i+1}$ are adjacent
for $1\le i<n$.  We write $\overline{K}$ for the equivalence class of
all isometries from $J$ into $S$ that are associate to $K$.

Remark. The term associate is closely related to the same term in
algebraic group theory, where two parabolic subgroups are called
associate if their Levi factors are conjugate. There is a small
difference because a standard parabolic subgroup corresponds to
a subdiagram of the Dynkin diagram, while we work with
{\it labeled} subdiagrams of a Coxeter diagram.

Suppose that $S$ is a spherical Coxeter diagram. Recall that $R$
is a normal subgroup of $\Gamma_J\subseteq Aut(J)$. We define an isometry
$K:J\mapsto S$ to be {\bf $R$-reflective} if $K$ is adjacent to $r(K)$
for some $r\in R$. We define the equivalence class $\overline{K}$ to
be {\bf $R$-reflective} if at least one of its elements is
$R$-reflective. (The reason for the name ``$R$-reflective'' appears in
lemmas 3.8 and 3.9.)

{\bf Example 2.3.} Suppose $J$ is $A_1$ and $S$ is $A_3$, with the
isometries from $J$ into $S$ labeled as $K_1$, $K_2$, $K_3$ in the
obvious way. Then the isometry $K_1$ is adjacent to $K_2$ as
$\sigma_{K_1(J)\cup K_2(J)}\sigma_{K_1}$ takes $K_1$ to
$K_2$. Similarly $K_2$ is adjacent to $K_3$, and $K_1$ is associate to
$K_3$ but not adjacent to $K_3$. The isometries $K_1$ and $K_3$ are
adjacent to themselves and $K_2$ is not adjacent to itself, so the
isometry $K_2$ is not $R$-reflective but the isometries $K_1$ and
$K_3$ are.  So $K_2$ is not $R$-reflective but the equivalence class
$\overline{K_2}$ is.

{\bf Example 2.4.} Suppose $K$ is $D_5$ and $S$ is $D_6$.
Then there
are exactly two isometries $K_1,K_2:J\mapsto S$ of $S$, which are adjacent
to each other but not to themselves. These two isometries are exchanged by
the nontrivial automorphism of $D_5$. Hence if $R$ contains just $1\in
Aut(D_5)$ then $K_1$, $K_2$, and the equivalence class $\overline
{K_1}=\{K_1,K_2\}$ are not $R$-reflective, but if $R$ is the whole of
$Aut(D_5)=Z/2Z$ then all of them are $R$-reflective.

{\bf Example 2.5} Suppose $K$ is $A_3$ and $S$ is $D_5$.
Then there
are 8 isometries $K:J\mapsto S$. These form two equivalence classes
under the relation of being associate, one of size 2 and one of size 6.
This shows that two isometries from
a connected diagram $J$ into $S$ need not be associate to conjugates of each
other under $Aut(J)$.

{\bf Example 2.6} Suppose $J$ is $A_2$. If $S$ is $A_n$ for $n\ge 2$
then there are two equivalence classes of isometries $K:J\mapsto S$, which
are exchanged by the nontrivial automorphism of $A_2$. However if $S$
is $D_n$ $(n\ge 4)$, $E_6$, $E_7$, or $E_8$ then there is only one
equivalence class, as the $A_2$ can be reversed by doing a ``three
point turn'' around the point of valence 3 in $S$.

We define a poset $P^+_3$ as follows. The objects of $P^+_3$
are pairs $(S,\overline{K})$ consisting of a spherical
subdiagram
$S$ of $\Pi$ and an equivalence class $\overline{K}$
of isometries from $J$ into $S$.
We define the partial order on $P^+_3$ by putting
$(S,\overline{K})\le (S',\overline{K'})$ if $S\subseteq S'$ and
$\overline{K}\subseteq \overline{K'}$. We define $P_3$ to be the
sub-poset of $P^+_3$ of elements $(S,\overline{K})$ such that
$\overline{K}$ is not $R$-reflective.

Note that the condition that $\overline{K}$ is not $R$-reflective
is quite restrictive and implies that $S$ is usually not much larger than 
$J$ and in any case has at most twice the rank of $J$. In particular
$S$ cannot contain any root orthogonal to $K(J)$ as this implies
that $K$ is $R$-reflective for any $R$. 

Suppose that $P$ is a poset acted on by a group $G$.
We define the {\bf homotopy quotient} of
$P$ by $G$ to be the following category $Q$. The objects of $Q$
are the elements of $P$. The morphisms from $p_1\in P$ to $p_2\in P$
correspond to the group elements $g\in G$ such that $g(p_1)\le p_2$,
and composition of morphisms is given by multiplication of group elements.
If $G$ is trivial this is the usual category associated to the poset
$P$, and if $P$ has just one point this is the usual category with one object
associated to the group $G$. If we take a full set of representatives
of the orbits of $G$ on $P$, then the full sub category of $Q$
with these objects is a skeleton of the category $Q$.

The poset $P_3$ is acted on by the group
$\Gamma_J\times \Gamma_\Pi \subseteq Aut(J)\times Aut(\Pi)$.
We define $Q_3$ to be
the homotopy quotient of $P_3$ by $\Gamma_J\times \Gamma_\Pi$, and
we construct the category $Q_4$ as a skeleton of $Q_3$ as above. In
other words the objects of $Q_4$ are a complete set of representatives
for the orbits of $\Gamma_J\times
\Gamma_\Pi$ on the elements of $P_3$ and the morphisms from
$(S,\overline{K})$ to $(S',\overline{K'})$ correspond to group elements
$\gamma\in \Gamma_J\times \Gamma_\Pi$ 
such that $\gamma((S,\overline{K}))\le (S',\overline{K'})$.

The main result of this paper is the following description of
the classifying space of $\Gamma_{\Omega}$.

\proclaim Theorem 2.7.
Suppose we are given the following objects.
\item{$(W_\Pi, \Pi)$} A Coxeter system.
\item{$\Gamma_\Pi$} A subgroup of $Aut(\Pi)$.
\item{$(W_J, J)$} A spherical Coxeter system with $J\subseteq \Pi$.
\item{$\Gamma_J$} A subgroup of $Aut(J)$.
\item{$R$} A normal subgroup of
$\Gamma_J$.

{\sl
Define $W_\Omega$, $\Gamma_\Omega$, and $Q_4$ as above, so that
the group $N_{W_\Pi.\Gamma_\Pi}(W_J;\Gamma_J)$ has the structure
$$W_J. W_\Omega.\Gamma_\Omega$$ where $W_\Omega$ is a Coxeter
group. Then  the component of $Q_4$ containing the object
$(J,\overline{id_J})$ is a classifying category for the group
$\Gamma_{\Omega}$.}

Theorem 2.7 gives a presentation of the group $\Gamma_{\Omega}$
because $\Gamma_\Omega$ is the fundamental group
of the category $Q_4$ with respect to the
basepoint $(J,\overline{id_J})$, and
it is easy to write down a presentation of the fundamental
group of any connected category $Q$ as follows. Choose a spanning
tree $T$ for the underlying 1-complex of $Q$ (which has a point for
each object of $Q$ and a 1-cell for each morphism). Then the
fundamental group $\Gamma_\Omega$ of $Q$ has a presentation as
follows. The group $\Gamma_\Omega$ has a generator $\overline{ g}$ for each
morphism $g$ of $Q$. The relations are $\overline{gh}=\overline
{g}\overline{h}$ whenever $gh$ is defined, and $g=1$ for $g$ in the
spanning tree.

{\bf Example 2.8.} Suppose we put $R=\Gamma_J=Aut(J)$ and
$\Gamma_\Pi=1$. Then we see from theorem 2.7 that
$N_{W_\Pi}(W_J)=W_J.W_{\Omega}.\Gamma_{\Omega}$, where
$\Gamma_{\Omega}$ is the fundamental group of the component of $Q_4$
corresponding to $J$. So in particular theorem 2.7 describes
normalizers of finite parabolic subgroups of Coxeter groups.  More
generally, Brink and Howlett [B-H] have described a presentation for
normalizers of possibly infinite parabolic subgroups of Coxeter
groups.  It is not trivial to see that the presentation given by [B-H]
is equivalent to the one given by theorem 2.7 (though of course this
follows from the fact that they are both presentations of the same
group).  Howlett pointed out to me that their result only requires
considering subdiagrams $S$ of $\Pi$ whose rank is at most $2+rank(J)$
to get the relations of $\Gamma_\Omega$, and of rank at most
$1+rank(J)$ to get the generators of $\Gamma_\Omega$.  It seems
possible that a similar simplification could be made to theorem 2.7 if
all that is required is a presentation rather than a classifying
space.  Perhaps the natural map from $\pi_i(Q_4^j)$ to
$\pi_i(Q_4)$ is an isomorphism for $i<j$ and an epimorphism for $i=j$,
where $Q_4^j$ is the full sub-category of $Q_4$ whose objects are the
elements $(S,\overline K)$ such that $rank(S)\le rank(J)+j$.  If so,
the map from $\pi_1(Q_4^2)$ to $\pi_1(Q_4)$ would be an
isomorphism, so this would give a closer connection to the
presentation of Brink and Howlett.  Their result also suggests that
theorem 2.7 could be generalized by allowing $W_J$ to be infinite and
modifying the definition of $Q_4$ to allow subdiagrams $S$ such that
$W_J$ has finite index in $W_S$.

\proclaim
3.~Proof of main theorem.

This section gives the proof of the 
theorem 2.7. The idea of the proof is to construct categories and 
functors according to the following diagram.
$$Q_1\longleftarrow Q_2\longrightarrow Q_3 \longleftarrow Q_4$$ It is
easy to show that a component of $Q_1$ is a classifying category for
$\Gamma_\Omega$. We also show that the functors between the categories
are all homotopy equivalences, so a component of $Q_4$ is a
classifying category for $\Gamma_\Omega$, which is what we wanted to
prove.

We define an isometry from $J$ into the roots of $W_\Pi$ to 
be {\bf primitive} if it is conjugate under $W_\Pi$ to
an isometry from $J$ into $\Pi$. An example of a non-primitive 
isometry is an isometry from $A_1^4$ into the roots of $D_4$.

We define a category $Q_1$ as follows. We define the poset $P_1$ to
be the poset of pairs $(C,K)$ where $K$ is a primitive
isometry from $J$ into the (possibly non-simple)
roots of $W_\Pi$, and $C$ is a Weyl chamber of
the reflection group $W_{\Omega_K}$ of $V_\Pi^{*W_K}$.
The partial order on $P_1$ is the trivial one with $(C_1,K_1)\le (C_2,K_2)$
if and
only if $(C_1,K_1)=(C_2,K_2)$.
The objects of $P_1$ are acted on in the natural
way by the group $W_\Pi.\Gamma_\Pi$ via its action on $V_\Pi$, and by the
group $\Gamma_J$ via its action on $J$. We define the category $Q_1$
to be the homotopy quotient of $P_1$ by the group
$\Gamma_J\times W_\Pi.\Gamma_\Pi$.

\proclaim Lemma 3.1.
The component of $Q_1$ containing the object $id_J:J\mapsto J\subseteq \Pi$
is a classifying category for the group
$\Gamma_{\Omega}$.

Proof. This follows immediately from the fact that $Q_1$ is a groupoid
such that the automorphism group of the object $K$
is the group $\Gamma_{\Omega_K}$.
This proves lemma 3.1.

Let $C_\Pi$ be the Weyl chamber of $W_\Pi$ defined in section 2. 
By a face of $C_\Pi$ we mean a nonempty intersection of
$C_\Pi$ with some of the hyperplanes bounding $C_\Pi$. The faces
of $C_\Pi$ of codimension $n$ correspond to the
spherical subdiagrams of $\Pi$ of rank $n$. We define a $\Pi$-cell 
to be a conjugate of a face of $C_\Pi$ under $W_\Pi$. The cone $X$ is
the union of all $\Pi$-cells, and the intersection of two $\Pi$-cells
is either empty or another $\Pi$-cell.

We define a category $Q_2$ and posets $P_2$, $P_2^+$ as follows.  The
objects of the poset $P^+_2$ are the pairs $(D,K)$ where $K$ is a
primitive isometry from $J$ into the roots of $W_\Pi$, and $D$ is a
$\Pi$-cell contained in $V_\Pi^{*W_K}$.  We define the partial order
on $P^+_2$ by saying $(D_1, K_1)\le (D_2, K_2)$ if $D_{1}\subseteq
D_{2}$ and $K_1=K_2$. We define $P_2$ to be the sub-poset of $P^+_2$
of elements $(D,K)$ such that $D$ is not contained in a reflection
hyperplane of $W_{\Omega_K}$.  The category $Q_2$ is defined to be the
homotopy quotient of $P_2$ by $\Gamma_J\times W_\Pi.\Gamma_\Pi$.

\proclaim Lemma 3.2. Suppose $G$ is a group, $P_1$ and $P_2$ 
are $G$-posets,  and $f$ is a morphism
of $G$-posets from  $P_2$ to
$P_1$. Also
suppose that for any $Y\in P_1$ the poset $f^{-1}(Y)$ is contractible
(in other words the corresponding simplicial complex is contractible).
Then the functor induced by $f$ between the homotopy quotient
categories $Q_2$, $Q_1$ of the posets
$P_2$ and $P_1$ by the group $G$ is a homotopy equivalence.

Proof. If $f$ is a functor from a category $Q_2$ to a category $Q_1$
and $Y$ is an object of $Q_1$ then we write $f^{-1}(Y)$ for the fiber
of $f$ over $Y$, in other words the sub category of $Q_2$ whose
morphisms are those mapped to the identity of $Y$ by $f$. We write
$Y\backslash f$ for the category consisting of pairs $(X,v)$ with
$v:Y\mapsto f(X)$, where a morphism from $(X,v)$ to $(X',v')$ is a
morphism $w:X\mapsto X'$ such that $f(w)v=v'$. Then a result due to
Quillen (the corollary to theorem A on page 9 of [Q]) states that $f$
is a homotopy equivalence provided that for all $Y$ in $Q_1$
the poset $f^{-1}(Y)$ is contractible
and the functor from $f^{-1}(Y)$ to $Y\backslash f$
taking $X$ to $(X,id_Y)$ has a right adjoint. (Here $id_Y$
is the identity morphism of $Y$.)

We will use Quillen's result to show that $f$ is a homotopy
equivalence. For any object $Y$ of $P_1$ the category $f^{-1}(Y)$ is
just the category of the poset $f^{-1}(Y)$, which is contractible by
assumption. So it only remains to check the condition about the existence
of a right adjoint from $Y\backslash f$ to $f^{-1}(Y)$. The category
$Y\backslash f$ has as objects pairs $(X,v)$ with $v\in G$, $v(Y)\le
f(X)$ and there is a morphism from $(X,v)$ to $(X',v')$ if and only if
$v^{-1}(X)=v'^{-1}(X')$, in which case the morphism is unique. We
define a functor $g$ from $Y\backslash f$ to $f^{-1}(Y)$ on objects by
$g((X,v))= v^{-1}(X)$. It is easy to check that this extends in a
unique way to morphisms. It is a right adjoint to $f$ because $Y\le
g((X,v))$ if and only if there is a morphism (necessarily unique) from
$f(Y)$ to $(X,v)$, both conditions being equivalent to $v(Y)\le
X$. This shows that the conditions of Quillen's result are satisfied,
so $f$ is a homotopy equivalence. This proves lemma 3.2.

\proclaim Lemma 3.3. The functor $f$ is a homotopy equivalence
from $Q_2$ to $Q_1$.

Proof. By lemma 3.2 it is sufficient to check that for
each $Y\in P_1$, the sub poset $f^{-1}(Y)$ of $P_2$ is contractible.
The poset
$f^{-1}(Y)$ is the poset of a cell decomposition of
a convex cone in a real vector space. As any convex set is contractible,
the poset $f^{-1}(Y)$ is also contractible. This proves lemma 3.3.

\proclaim Lemma 3.4. Suppose $(W,S)$ is a spherical Coxeter system
acting on the vector space $V_S$ with Weyl chamber $C$. Suppose $K$ is
an isometry from $J$ into $S$. Let $V_S^{*W_K}$ be the subspace of
$V^*$ fixed by $W_K$, where $W_K$ is the reflection group
whose simple roots are the points $K(J)$. The walls of $C\cap
V_S^{*W_K}$ correspond to the points in $S$ not in the image of $K$;
let $s$ be one of these points and let $s^\perp\cap V_S^{*W_K}$ be the
wall in $V_S^{*W_K}$ corresponding to $s$. Choose $w\in W$
so that  $w(C)$ is
the (unique) Weyl chamber of $W$ such that $w(C)\cap V_S^{*W_K}$ is
the cell in $V_S^{*W_K}$ on the other side of $s^\perp\cap V_S^{*W_K}$
to $C\cap V_S^{*W_K}$ and such that $C$ and $w(C)$ are both in the
same Weyl chamber of $W_K$. Then
$$w=\sigma_{K(J)\cup s}\sigma_{K(J)}.$$

Proof. We can reduce to the case when $S=K(J)\cup s$, so that
$V_S^{*W_K}$ is one dimensional and $s^\perp\cap V_S^{*W_K}$ is
just the point 0. Then $\sigma_S(\sigma_{K(J)}(C))=-\sigma_{K(J)}(C)$ which
contains $-C\cap V_S^{*W_K}$, so $\sigma_S(\sigma_{K(J)}(C))\cap
V_S^{*W_K}$ is a cell on the other side of $s^\perp\cap V_S^{*W_K}$ to
$C\cap V_S^{*W_K}$. Moreover $\sigma_{K(J)}(C)$ is in the opposite
Weyl chamber of $K(J)$ to $C$, and $\sigma_S(\sigma_{K(J)}(C))$ is in
the opposite Weyl chamber to $\sigma_{K(J)}(C)$, so
$\sigma_S(\sigma_{K(J)}(C))$ is in the same Weyl chamber of $K(J)$ as
$C$. This shows that the element $w$ of the lemma is
$\sigma_S\sigma_{K(J)}$.
This proves lemma 3.4.

A result similar to the following lemma 
(using subsets of $S$ rather than isometries $K:J\mapsto S$) is given in [H,
lemma 5] when $S$ is finite and in [D] for arbitrary $S$.

\proclaim Lemma 3.5. Suppose $(W,S)$ is a Coxeter system,
$J$ is a spherical Coxeter diagram, and
$K$ and $K'$ are two isometries from $J$
into $S$. Then $K$ and $K'$ 
are conjugate under
$W$ if and only if they are associate.

Proof. 
First suppose that $K$ and $K'$ are adjacent. Then
$\sigma_{K\cup K'}(\sigma_{K}(K))=K'$, so $K$ and $K'$ are conjugate
under $W$. Next suppose $K$ and $K'$ are associate. Then by
definition we can find a sequence $K=K_1$, $K_2,\ldots, K_n=K'$ such
that $K_i$ and $K_{i+1}$ are adjacent for all $i$. Hence $K=K_1$ and
$K'=K_n$ are also conjugate under $W$. So associate isometries from $J$ into
$S$ are conjugate under $W$.

Conversely, suppose that $K$ and $K'$ are conjugate by an element $w\in
W$. Consider the subspace $V_S^{*W_K}$ of
$V_S^*$, and the codimension 0 cells in it of the form
$V_S^{*W_K}\cap C$ for some Weyl chamber $C$ of $W$. For
any two such cells, for example $D=V_S^{*W_K}\cap C$ and
$D'=V_S^{*W_K}\cap w(C)$, we can find a sequence
$D=D_1,D_2,\ldots, D_n=D'$ such that $D_i$ and $D_{i+1}$ are adjacent
by a face of codimension 1 in $V_S^{*W_K}$. For each $i$
let $C_i$ be the (unique) Weyl chamber whose intersection with
$V_S^{*W_K}$ is $D_i$ and that is contained in the Weyl
chamber of $W_K$. We identify each $C_i$ with $C_\Pi$ using $w_i$.
The set $K(J)$ is a subset of the simple roots of $C_i$, so
$K_i=w_i^{-1}(K)$ maps $J$ to the simple roots of $C_\Pi$. The isometries
$K_i$ and $K_{i+1}$ are adjacent for all $i$, because the element
$w_{i+1}w_i^{-1}$ mapping $C_i$ to $C_{i+1}$ is equal to
$\sigma_{K(J)\cup s_i}\sigma_{K(J)}$, where $s_i$ is the simple root
of $C_i$ orthogonal to $D_i\cap D_{i+1}$ but not to $D_i$. Therefore
$K=K_1$ and $K'=K_n$ are adjacent. This proves lemma 3.5.

Suppose that $P$ is a $W$-poset for a group $W$ with the property
that if $p\le w(p)$ for $w\in W, p\in P$ then $p=w(p)$. We define the 
{\bf quotient}
$W\backslash P$ of $P$ by $W$ to be the poset whose elements are
the orbits $Wp$ of $W$ acting on $P$, where we put $Wp\le Wq$ if
$w(p)\le q$ for some $w\in W$.
This should not be confused with the homotopy quotient of $P$ by $W$.

\proclaim Lemma 3.6. The $\Gamma_J\times \Gamma_\Pi$ posets
$P_3^+$ and $W_\Pi\backslash P_2^+$ are isomorphic.

Proof. We will construct an isomorphism $f$ of posets from
$W_\Pi\backslash P^+_2$ to $P^+_3$. Suppose $(D,K)$ is an element of
$P^+_2$ representing an element of $W_\Pi\backslash P^+_2$. We can
find an element $w$ of $W_\Pi$ such that $w(D)\subseteq C_\Pi$ and
$w(K(J))\subseteq \Pi$. We define $f((D, K))$
to be $(S,\overline{w(K)})$, where $S$ is the set of simple roots of
$C_\Pi$ orthogonal to $w(D)$. We check that this is well defined even
though $w$ is not unique. To prove this we can assume that $D\subseteq
C_\Pi$ and $K(J)\subseteq \Pi$. Then the different possibilities for
$w$ are elements of the group generated by the reflections fixing $D$
and the Weyl chamber of $K(J)$. But these elements take $K$ to an
associated isometry $K:J\mapsto S$, so the equivalence class $\overline{K}$
is well defined by lemma 3.5. These elements also take $S$ to $S$, so
$S$ is well defined. This proves that $(S,\overline{K})$ is uniquely
defined.

The isomorphism $f$ of posets from $W_\Pi\backslash P^+_2$
to $P^+_3$ obviously preserves the $\Gamma_J\times \Gamma_\Pi$ action on
both posets. This proves lemma 3.6.

\proclaim Lemma 3.7. Suppose $(W,S)$ is a spherical Coxeter system,
$K$ is an isometry $K:J\mapsto S$, and $s$ is a point of $S$ not in
$K(J)$. Then there is an element of $W$ mapping $V_S^{*W_K}$ to itself
and acting on $V_S^{*W_K}$ as reflection in $s^\perp\cap V_S^{*W_K} $
if and only if $-\sigma_{K(J)\cup s}$ maps $K(J)$ to itself.  If such
an element of $W$ exists, then there is a unique such element $w$
mapping $K$ to itself, given by $w=\sigma_{K(J)\cup s}\sigma_{K(J)}$.

Proof. If an element of $w$ maps $V_S^{*W_K}$ to itself then there is
a unique element of $W$ with the same action on $V_S^{*W_K}$ and
mapping $K(J)$ to itself because $W_K$ acts simply transitively on its
Weyl chambers, so we may assume that $w$ maps $K(J)$ to itself.  If in
addition $w$ acts on $V_S^{*W_K}$ as reflection in $s^\perp\cap
V_S^{*W_K} $ then by lemma 3.4 $w$ must be $\sigma_{K(J)\cup
s}\sigma_{K(J)}$.

Conversely if $-\sigma_{K(J)\cup s}$ maps $K(J)$ to itself then
$\sigma_{K(J)\cup s}\sigma_{K(J)}$ maps $K(J)$ to itself and acts on
$V_S^{*W_K}$ as reflection in $s^\perp\cap V_S^{*W_K} $.  This proves
lemma 3.7.

\proclaim Lemma 3.8. Suppose $K$ is an isometry from $J$ into a spherical
subdiagram $S$ of $\Pi$. Then $K:J\mapsto S$ 
is $R$-reflective if and only if the
$\Pi$-cell $S^\perp\cap C_\Pi$ is contained in a reflection hyperplane
of $W_{\Omega_K}$ of the form $s^\perp\cap V_\Pi^{*W_K}$ for $s\in \Pi$. 

Proof. First suppose that $K:J\mapsto S$ is $R$-reflective. 
Then there is a point $s\in S$ not in $K(J)$ such that
$w=\sigma_{K(J)\cup s}\sigma_{K(J)}$ acts on $K(J)$ as an element of $R$. 
So $w$ is a reflection of $W_\Omega$ corresponding to the hyperplane
$s^\perp\cap V_\Pi^{*W_K}$, and this hyperplane contains  $S^\perp\cap C_\Pi$.

Conversely, suppose that $S^\perp\cap C_\Pi$ is contained in a
reflection hyperplane of $w\in W_{\Omega_K}$ of the form $s^\perp\cap
V_\Pi^{*W_K}$ for $s\in \Pi$. Then we must have $s\in S$ because
$S^\perp\cap C_\Pi\subseteq s^\perp$. The element $w$ must be equal to 
$\sigma_{K(J)\cup s}\sigma_{K(J)}$, and this element acts on 
$K$ as an element of $R$ because $w\in W_\Omega$. Therefore
$K:J\mapsto S$ is $R$-reflective. 
This proves lemma 3.8.

\proclaim Lemma 3.9. Suppose $K$ is an isometry from $J$ into a spherical
subdiagram $S$ of $\Pi$. Then the equivalence class $\overline{K}$ is
$R$-reflective if and only if the $\Pi$-cell $S^\perp\cap C_\Pi$ of
$V_\Pi^{*W_K}$ corresponding to $S$ is contained in a reflection
hyperplane of $W_{\Omega_K}$.

Proof. Suppose the cell of $V_\Pi^{*W_K}$ corresponding to $S$ is
contained in a reflection hyperplane of $W_{\Omega_K}$. Choose a Weyl
chamber $C_\Pi$ for $W_\Pi$ such that this reflection hyperplane is a
wall of $C_\Pi\cap V_\Pi^{*W_K}$. Then by lemma 3.8 the corresponding
isometry $K':J\mapsto S$ is $R$-reflective, and by lemma 3.5 is
associate to $K$. So $\overline K$ is $R$-reflective.

Conversely suppose that $\overline K$ is $R$-reflective.  Then $w(K)$
is $R$-reflective and has image in $S$ for some $w\in W_\Omega$.  By
lemma 3.8 this implies that $w(S^\perp\cap C_\Pi)$ is contained in a
reflection hyperplane of $W_\Omega$, so the same is true of
$S^\perp\cap C_\Pi$.  This proves lemma 3.9.

\proclaim Lemma 3.10. Suppose $P$ is a $G$-poset for some group $G$.
Suppose that $W$ is a normal subgroup of $G$
such that if $p\le q$ and $w(p)\le q$ for some $w\in W$,
$p,q\in P$, then $w=1$.
Then the homotopy quotient $Q_2$ of $P$ by $G$
is equivalent to the homotopy quotient
$Q_3$ of
$W\backslash P$ by $G/W$.

Proof.
Recall from [M, theorem 1, page 91]
that if $f$ is any functor from a category
$Q_2$ to a category $Q_3$ then $f$ is an equivalence if
the following two conditions are satisfied:
\item{1} Any object of $Q_3$ is isomorphic
to some object in the image of $f$.
\item{2} For any two
objects $p$, $q$ of $Q_2$, $f$ induces an isomorphism from
$Mor(p,q)$ to $Mor(f(p),f(q))$.

We will apply this to show that our categories $Q_2$ and $Q_3$ are
equivalent.  We define $f$ on objects by $f(p)=Wp$, and define $f$ on
morphisms using the obvious homomorphism from $G$ to $G/W$. Condition
1 above is satisfied because every element of $P_3$ is the image of an
element of $P_2$, so every object of $Q_3$ is the image of an object
of $Q_2$.  Suppose $p$ and $q$ are objects of $P_2$. The set of
morphisms of $Q_2$ from $p$ to $q$ can be identified with the set of
elements $g$ of $G$ such that $g(p)\le q$, and $Mor_{Q_3}(f(p),f(q))$
can also be identified with the set of elements $g$ of $G$ such that
$g(p)\le q$, so condition 2 above is satisfied. This shows that $f$ is
an equivalence and proves lemma 3.10.

\proclaim Lemma 3.11. If $p\le q$ and $w(p)\le q$ for some $w\in W_\Pi$,
$p,q\in P_2$, then $w=1$.

Proof. Suppose that $p=(D, K)$. Then $q=(D_1, K)$
for some $D_1$ containing $D$ as $p\le q$. But then $w(D)\subseteq D_1$,
so $w(D)=D$ as no two distinct subsets $D$, $w(D)$ of $D_1$ are conjugate
under $W_\Pi$, as $D_1$ is contained in a fundamental domain of $W_\Pi$.
Hence we can assume that $w$ fixes $D$ as well as $K$.

The subgroup of $W_\Pi$ fixing $D$ is a finite reflection group $W_D$
generated by the reflections of $W_\Pi$ fixing $D$ because $D\subseteq C_\Pi$. 
The subgroup of
$W_\Pi$ fixing $K$ is generated by the reflections fixing all elements
of $K(J)$. Any such reflection is in $W_{\Omega_K}$ because
$R$ contains 1. However, the condition that $p=(D,K)\in P_2$ implies
that  there are no reflections of $W_{\Omega_K}$ fixing $D$.
Hence the subgroup of $W_D$ fixing $K$ is trivial.
So any element $w\in W_\Pi$ such that $w(p)\le q$ is trivial.
This proves lemma 3.11.

\proclaim Lemma 3.12. The categories
$Q_2$ and $Q_3$ are equivalent.

Proof. Lemmas 3.10 and 3.11 show that there is an equivalence of
categories from $Q_2$ to the homotopy quotient of  $W_\Pi\backslash P_2$
by  $\Gamma_J\times\Gamma_\Pi$.
Lemma 3.6 shows that the $\Gamma_J\times \Gamma_\Pi$ posets $P_3^+$
and $W_\Pi\backslash P_2^+$ are isomorphic. Lemma 3.9 shows that the
subset $P_2$ of $P^+_2$ corresponds under this isomorphism to the
subset $W_\Pi\backslash P_2$ of $W_\Pi\backslash P^+_2$, so the
$\Gamma_J\times \Gamma_\Pi$ posets $P_3$ and $W_\Pi\backslash P_2$ are
isomorphic.  Therefore the category $Q_2$ is equivalent to the
homotopy quotient of $P_3$ by $\Gamma_J\times
\Gamma_\Pi$, which is just $Q_3$.
This proves lemma 3.12.

\proclaim Lemma 3.13. The natural injection from $Q_4$ to $Q_3$ is
an equivalence of categories.

Proof. This follows from the fact that $Q_4$ is a skeleton of $Q_3$,
so the natural injection is an equivalence of categories.
This proves lemma 3.13.

We can now prove theorem 2.7. By lemmas 3.3, 3.12, and 3.13,
the categories $Q_4$ and $Q_1$ are homotopy equivalent.
So by lemma 3.1, the component of $Q_4$ containing $(J, \overline{id_J})$
is a classifying category for $\Gamma_{\Omega}$. This
proves theorem 2.7.

\proclaim 4.~The structure of $\Gamma_{\Omega}$.

\proclaim Theorem 4.1. The kernel of the natural map from
$\Gamma_{\Omega}$ to $\Gamma_J\times \Gamma_\Pi$ has finite
cohomological dimension.

Proof. The classifying category of the kernel is just a component of
the category of the poset $P_3$.
A case by case check on possible Coxeter diagrams 
shows that the lengths of chains in $P_3$ are
bounded (by $rank(J)+1$ for example), so the corresponding simplicial
complex has finite dimension at most $rank(J)$. Therefore the kernel 
has cohomological dimension  at most $rank(J)$. 
This proves
theorem 4.1.

\proclaim Corollary 4.2. If $\Gamma_\Pi$ has finite virtual cohomological
dimension, then so does $\Gamma_{\Omega}$.

Proof. This follows immediately from theorem 4.1 and the fact that
$\Gamma_J$ is finite and standard properties of the virtual
cohomological dimension.

Howlett showed that if $W_\Pi$ is finite then the group
$\Gamma_{\Omega}$ is a subgroup of $Aut(J)\times Aut(\Pi)$. We can
deduce this from theorem 4.1 as follows. If $W_\Pi$ is finite then so
is the kernel of the map from $\Gamma_{\Omega}$ to $Aut(J)\times
Aut(\Pi)$. On the other hand this kernel has finite cohomological
dimension by theorem 4.1. But any finite group of finite cohomological
dimension must be trivial, so the natural map from $\Gamma_{\Omega}$
to $Aut(J)\times Aut(\Pi)$ is injective. If $W_\Pi$ is infinite then
this kernel is usually infinite, as can be seen from most of the
examples below. The fact that this kernel no longer vanishes is the
main reason why normalizers of parabolic subgroups of Coxeter groups
are more complicated to describe when the Coxeter group is infinite.

\proclaim 5.~Examples.

{\bf Example 5.1.} Suppose that $J$ is $A_1$ and the group
$\Gamma_\Pi$ is trivial.  In this case Brink [Br] gave an elegant
description of the group $\Gamma_{\Omega}$ as follows. Form the graph
obtained from the Coxeter graph $\Pi$ by keeping only the edges of odd
order. Then for any point $J=A_1$ of this new graph, the centralizer
of the corresponding reflection (which is the group $\Gamma_{\Omega}$
corresponding to $J$)  is the fundamental group of this graph with
basepoint the chosen point, and in particular $\Gamma_{\Omega}$ is a
free group.

We now check that this is equivalent to the description given by
theorem 2.7. The only subdiagrams $S$ with a non $R$-reflective class
$\overline{K}$ are the points of the Coxeter graph or the edges of odd
order together with their endpoints. So the category $Q_4$ has an
object for each point or odd order edge of the Coxeter graph. The only
non-identity morphisms correspond to inclusions of points in
edges. The classifying space of this category is just the first
barycentric subdivision of Brink's graph. In particular the
fundamental group of this category with some object as basepoint is
canonically isomorphic to the fundamental group of Brink's graph with
the corresponding point as basepoint. This verifies that the
description of $\Gamma_{\Omega}$ in [Br] is equivalent to
the description in theorem 2.7.

The following lemma can often be used to find the fundamental group
of a category with at most 2 objects. 

\proclaim Lemma 5.2. Suppose $A$ and $B$ are subgroups of a group.
Let $Q$ be the category with 2 objects $p$ and $q$ such that
$Mor(p,p)=A$, $Mor(p,q)=BA$, $Mor(q,q)=B$, $Mor(q,p)=\emptyset$,
with composition defined in the obvious way. Then
$\pi_1(Q)=A*_{A\cap B} B$.

Proof. This can be proved by writing down a set of generators and relations
for the fundamental group, and checking that they are equivalent
to a set of generators and relations for $A*_{A\cap B} B$. We will leave
the details to the reader.

For most of the examples below we will take $W_\Pi.\Gamma_\Pi$ to be
the group of automorphisms of the even 26 dimensional Lorentzian
lattice $II_{1,25}$ not exchanging the two cones of norm 0 vectors.
According to Conway [C], the Coxeter group $W_\Pi$ has a simple
reflection $r_\lambda\in \Pi$ for each vector $\lambda$ of the Leech
lattice $\Lambda=\Pi$, and the order of $r_\lambda r_\mu$ is 1, 2, 3,
or $\infty$ according to whether $(\lambda-\mu)^2$ is 0, 4, 6, or
greater than 6. The group $\Gamma_\Pi$ is the automorphism group
$\Lambda.Aut(\Lambda)=\Lambda.(Z/2Z).Co_1$ of the affine Leech
lattice, where $\Lambda$ is the subgroup of translations and $Co_1$ is
Conway's largest sporadic simple group.

If $J$ is a spherical subdiagram of $\Lambda$ then there is a
homomorphism from $N_{W_\Pi.\Gamma_\Pi}(W_J)/W_J$ to the automorphism
group of the lattice $J^\perp$. This has finite kernel and co-finite
image, so theorem 2.7 can usually be used to describe the automorphism
group of the lattice $J^\perp$. Most of the remaining examples in this
section use this idea.

{\bf Example 5.3.} Suppose $L$ is the even Lorentzian lattice of dimension
20 and determinant 3. Let $W^{(2)}(L)$ be the subgroup of $Aut(L)$ generated
by reflections of norm $-2$ vectors of $L$.
Vinberg showed in [V] that $Aut(L)^+/W^{(2)}(L)$
was the automorphism group of a certain K3 surface modulo a cyclic
subgroup, and
also showed that this group was an extension of a group of order 72 by
a free product of 12 groups of order 2. We will show how to recover
Vinberg's description of $Aut(L)^+/W^{(2)}(L)$ from theorem 2.7.

We take $J$ to be an $E_6\subset \Lambda$ so that $L=J^\perp$, and
take $R=1\subset
\Gamma_J=Aut(E_6)=Z/2Z$.  The category $Q_4$ contains exactly two
objects, corresponding to an $E_6$ and an $E_7$ in $\Lambda$.  For
each subdiagram $X$ of $\Lambda$ we write $G(X)$ for the automorphisms
of $\Pi=\Lambda$ mapping $X$ into itself.  The morphisms from the
$E_6$ object to itself form a group $G(E_6)$ of order 72. The
morphisms from the $E_7$ object to itself form a group $Z/2Z\times
G(E_7)$ of order $2\times 6=12$ (where the $Z/2Z$ comes from the group
$\Gamma_J$). The morphisms from $E_6$ to $E_7$ can be identified with
$Z/2Z\times G(E_6)$. By lemma 5.2 the group $\Gamma_{\Omega}$ is
isomorphic to $G(E_6)*_{G(E_7)}(G(E_7)\times Z/2Z)$. If $A$, $B$, and
$C$ are any groups with $B\subseteq A$ then $A*_B(B\times C)$ is a
semidirect product of a normal subgroup isomorphic to the free product
of $|A|/|B|$ copies of $C$ and with the quotient by this normal
subgroup isomorphic to $A$. Hence $\Gamma_{\Omega}$ has a normal
subgroup isomorphic to the free product of 12 copies of $Z/2Z$, and
the quotient is the group $G(E_6)$ of order $72$. This is equivalent
to the description of this group given by Vinberg in [V]. In fact we
can describe the various parts of Vinberg's description in terms of
the Leech lattice as follows: the group of order 2 is the group of
automorphisms of $E_6$, the group of order 72 is the subgroup of
$Aut(\Lambda)$ mapping an $E_6$ into itself, and the number 12 is the
number of $E_7$'s of $\Lambda$ containing an $E_6$.

The group $Aut(L)^+$ also contains reflections in norm $-6$ vectors.
The quotient by the full reflection group is finite of order 72,
isomorphic to $G(E_6)$. In this case the category $Q_4$ contains just
one point. Note that the reflections of norm $-6$ vectors induce the
nontrivial automorphism of $E_6$. The 12 elements of order 2 in the
paragraph above are in fact reflections of norm $-6$ vectors.  See [V]
or [B] for more details of this case.

If $L$ is the even  Lorentzian lattice 
of determinant 4 and dimension  20,
which is again the Picard lattice of a $K3$ surface,
then Vinberg gave a similar description of the automorphism  group
as an extension $((Z/2Z)*(Z/2Z)*(Z/2Z)*(Z/2Z)*(Z/2Z)).S_5$
of the symmetric group $S_5$ by the free product of
5 group of order 2. 
This group can also be calculated using theorem 2.7. 
The corresponding category $Q_4$ has 2 objects,
corresponding to a $D_6$ or $D_7$ in $\Lambda$, the group $S_5$ is
the subgroup of $Aut(\Lambda)$ mapping the $D_6$ to itself, the group 
$Z/2Z$ is the group of automorphisms of the $D_6$ diagram, and
the number 5 of copies of $Z/2Z$  is the number of $D_7$'s containing a $D_6$.

{\bf Example 5.4} As a more complicated example we will describe the
automorphism group of the Picard lattice $L$ of the ``next most
algebraic K3 surface'', in other words $L$ is the 20 dimensional even
Lorentzian lattice of determinant 7. We take $\Pi=\Lambda$,
$\Gamma_\Pi=\Lambda.Aut(\Lambda)$, $J=A_6$, $R=1$,
$\Gamma_J=Aut(A_6)=Z/2Z$. Then $\Gamma_{\Omega}$ is the subgroup of
elements of $Aut(L)$ fixing a Weyl chamber of the reflection group
generated by the reflections of norm $-2$ vectors. By theorem 2.7 the
group $\Gamma_{\Omega}$ is the fundamental group of the category
$Q_4$. The category $Q_4$ has exactly 5 objects, corresponding to the
5 orbits of Coxeter diagrams $A_6$, $A_7$, $D_7$, $E_7$, and $D_8$
with a non $R$-reflective isometry from $A_6$ into them. (The group
$\Lambda.Aut(\Lambda)$ acts transitively on the embeddings of any of
these Coxeter diagrams into $\Lambda$.) Note that for $D_8$ and $E_7$
there is only one equivalence class of isometries from $A_6=J$ into
it, while for $A_6$, $A_7$, and $D_7$ there are two classes, which are
exchanged by $Aut(J)=Z/2Z$. The category $Q_4$ looks like this.
$$\matrix{
&&&A_7(48)\cr
&&(672)\nearrow&&\searrow(192)\cr
E_7(12)&{\longleftarrow}
&A_6(336)&{\longrightarrow}&D_8(16)\cr
&(672)&(672)\searrow&(2016)&\nearrow(48)\cr
&&&D_7(24)\cr
}$$
Here the numbers are the numbers of morphisms between pairs of objects in
$Q_4$.

{\bf Example 5.5.} Kondo in [K] studied the automorphism group of a
generic Jacobian Kummer surface by embedding its Picard lattice $L$ as
the orthogonal complement of a certain $J=A_3A_1^6$ in $II_{1,25}$,
and used this to describe a generating set for the automorphism  group.
(Note that the Leech lattice contains more than 1 orbit of subdiagrams
of the form $A_3A_1^6$; the one used by Kondo has largest possible
stabilizer in $Aut(\Lambda)$.) By results of Nikulin [N] the
automorphism group of the K3 surface is the subgroup of
$Aut(L)^+/W^{(2)}(L)$ of elements acting on $L'/L$ as $\pm 1$.  This
is just the group $\Gamma_\Omega $ of theorem 2.7 where we take
$\Pi=\Lambda$, $\Gamma_\Pi=\Lambda.Aut(\Lambda)$, $J$ to be Kondo's
$A_3A_1^6$, $R=1$, and $\Gamma_J =$ the subgroup of order 2 of
$Aut(J)$ generated by the nontrivial automorphism of the $A_3$.  The
category $Q_4$ is not connected and the component containing $J$ seems
quite complicated. Some partial calculations I have done suggest that
the automorphism group of the generic Jacobian Kummer surface might be
$$(W.(Z/2Z)^5)*(Z/2Z)*(Z/2Z)*(Z/2Z)*(Z/2Z)*(Z/2Z)*(Z/2Z)$$
where $W$ is a Coxeter group of rank $32+60$ generated 
by the 16 projections, 16 correlations, sixty Cremona transformations
and the $(Z/2Z)^5$ is generated by sixteen translations and a switch [K], 
but I have not proved this
rigorously. 

It is much easier to work out the group $\Gamma_\Omega$ for $R$ and $\Gamma_J$
replaced by $Aut(J)=Z/2Z\times S_6$. In this case the group $\Gamma_\Omega$
has a normal subgroup of index $|S_6|$ isomorphic to the quotient of 
the automorphism  group of a generic Jacobian Kummer surface 
by the Coxeter group generated by reflections in norm $-4$ vectors. 
Theorem 2.7 describes the
$\Gamma_\Omega$ as the fundamental group of a component of the finite
category $Q_4$. This component
has just two elements,
corresponding to the Coxeter diagrams $A_3A_1^6$ and $A_5A_1^5$. Using
lemma 5.2 and the results in [K] we see that
$$\Gamma_\Omega=((Z/2Z)^5.S_6)*_{S_5}(S_5\times Z/2Z)$$ where
$(Z/2Z)^5.S_6$ is the subgroup of $\Gamma_\Pi$ fixing $A_3A_1^6$ ([K,
lemma 4.5]), and $S_5\times Z/2Z$ is the subgroup of $\Gamma_\Pi$
fixing $A_5A_1^5$.
In particular $\Gamma_\Omega$ has a normal subgroup which is the free
product of 192 groups of order 2, and the quotient by this normal
subgroup is $(Z/2Z)^5.S_6$.  If we change $R$ to $1$ but keep
$\Gamma_J=Aut(J)$ then $\Gamma_\Omega$ becomes the group
$Aut(L)^+/W^{(2)}(L)$ which appears to be
$$(W.(Z/2Z)^5.S_6)*_{S_5}(S_5\times (Z/2Z))$$
though I have not proved this rigorously. 

For more examples of automorphism  groups of Kummer surfaces,
corresponding to the cases $J=D_4^2$, $D_4A_3$, $D_4A_2$, or $A_3^2$, 
see [K-K].

{\bf Example 5.6} Suppose that $J$ in theorem 2.7 contains no
components of types $A_n$ $(n\ge 1)$ or $D_5$, and assume that
$R=\Gamma_J=Aut(J)$. Then the map from $\Gamma_{\Omega}$ to
$\Gamma_J\times \Gamma_\Pi$ is injective. This follows because a case
by case check over all irreducible spherical Coxeter diagrams shows
that any isometry from $J$ into a strictly larger spherical Coxeter
diagram is $R$-reflective.

{\bf Example 5.7} We show how to explain Vinberg's result [V, V-K] 
that the reflection group of $I_{1,n}$ has finite index if and only if
$n\le 19$. Following Conway and Sloane [C-S] we write the even
sublattice $L$ of $I_{1,n}$ as $D_{25-n}^\perp$ in $II_{1,25}$ for
$n\le 23$, where $D_3=A_3$ and $D_2=A_1^2$. In theorem 2.7 we take
$\Pi=\Lambda$, $\Gamma_\Pi=\Lambda.Aut(\Lambda)$, $J=D_{25-n}$,
$R=\Gamma_J$ a subgroup of order 2 of $Aut(J)$ (which is equal to
$Aut(J)$ for $n\ne 21$). Then the quotient of $Aut(L)^+$ by its
reflection subgroup is the group $\Gamma_{\Omega}$. For $n\le 19$ the
group $\Gamma_{\Omega}$ is finite by example 5.6. For $n=20$ this
argument breaks down because $J$ is the ``exceptional'' case $D_5$ of
example 5.6. For $20\le n\le 23$ we can still describe the group
$\Gamma_{\Omega}$ explicitly using theorem 2.7; see [B, theorem 6.6]
for details.  When $n=21$ this gives a natural example with
$\Gamma_J\ne Aut(J)$.

If we take $J$ to be $D_4$ and take $R=\Gamma_J$ to be the symmetric group
$S_3=Aut(D_4)$ instead of a group of order 2 then $
W_{\Omega}.\Gamma_\Omega$ is the automorphism group of the even
sublattice of $I_{1,21}$, and $\Gamma_{\Omega}$ is a finite group. See
[B, p. 149] for details.

{\bf Example 5.8.} The groups $\Gamma_{\Omega}$ have many of the
properties of arithmetic groups; for example, they often have finite
classifying categories. (It follows easily from [S] that arithmetic
groups have this property.) It is natural to ask when they are
arithmetic. There seems to be no obvious general way of deciding
this. The following argument can often be used to show that
$\Gamma_{\Omega}$ is not arithmetic. First of all a theorem due to
Margulis [Ma, page 3] implies that if a group is an arithmetic
subgroup of a group of rank at least 2, then all normal subgroups are
either in the center or of finite index. Secondly, a theorem of Borel
and Serre [B-S, 11.4.4] says that if a group $\Gamma$ is arithmetic in
a Lie group $G$ then $d=r+vcd(\Gamma)$ where $d$ is the dimension of
the symmetric space of $G$ and $r$ is the rank of $G$ and
$vcd(\Gamma)$ is the virtual cohomological dimension of $\Gamma$. Let
us use these results to prove that the group $\Gamma_{\Omega}$ of
example 5.3 (the automorphism group of a K3 surface) is not arithmetic
in any Lie group $G$. The group $G$ must have rank 1 by the theorem of
Margulis, as $\Gamma_{\Omega}$ has non abelian free subgroups of
finite index and so cannot be an arithmetic subgroup of a group of
rank at least 2. Its virtual cohomological dimension is one, so by the
theorem of Borel and Serre [B-S, 11.4.4] the symmetric space of $G$
must have dimension $1+1=2$. But any finite subgroup of
$\Gamma_{\Omega}$ must fix a point of this symmetric space, and
therefore acts faithfully on the 2 dimensional tangent space of this
point. But $\Gamma_{\Omega}$ has finite subgroups that are too large
to have 2 dimensional faithful representations. Hence
$\Gamma_{\Omega}$ is not arithmetic. Note that $\Gamma_{\Omega}$
has subgroups of finite index that are free and therefore arithmetic.

\proclaim References.

\item{[B]} R. E. Borcherds, Automorphism groups of Lorentzian lattices,
J.Alg. Vol 111 (1987) No. 1, pp. 133-153.
\item{[B-S]} A. Borel, J.-P. Serre, Corners and arithmetic groups.
Comment. Math. Helv. 48 (1973), 436--491.
\item{[B-H]} B. Brink, R. B. Howlett, 
Normalizers of Parabolic Subgroups in Coxeter Groups,
1998 preprint, available from http://www.maths.usyd.edu.au:8000/u/\~{}bobh/
\item{[Bo]} N. Bourbaki,
Groupes et alg\`ebres de Lie. Chapitres IV-VI.
Actualit\'es Scientifiques et Industrielles, No. 1337
Hermann, Paris 1968.
\item{[Br]} B. Brink, 
On centralizers of reflections in Coxeter groups. Bull. London Math. Soc. 28
(1996), no. 5, 465--470.
\item{[C]} J. H. Conway,
The automorphism group of the $26$-dimensional even unimodular Lorentzian
lattice. J. Algebra 80 (1983), no. 1, 159--163. 
Reprinted as chapter 27 of
``Sphere packings, lattices and groups.'' 
Grundlehren der Mathematischen Wissenschaften
290. Springer-Verlag, New York, 1993.  ISBN: 0-387-97912-3
\item{[C-S]} J. H. Conway, N. J. A. Sloane,
Leech roots and Vinberg groups. Proc. Roy. Soc. London Ser. A
384 (1982), no. 1787, 233--258.
Reprinted as chapter 28 of
``Sphere packings, lattices and groups.'' 
Grundlehren der Mathematischen Wissenschaften
290. Springer-Verlag, New York, 1993.  ISBN: 0-387-97912-3
\item{[D]} V. V. Deodhar, On the root system of a Coxeter group,
Comm. Algebra 10 (1982) 611--630.
\item{[Hi]} H. Hiller,
Geometry of Coxeter groups. Research Notes in Mathematics, 54. Pitman
(Advanced Publishing Program), Boston, Mass.-London, 1982. 
ISBN: 0-273-08517-4
\item{[H]} R. B. Howlett,
Normalizers of parabolic subgroups of reflection groups. J. London Math.
Soc. (2) 21 (1980), no. 1, 62--80.
\item{[Ke]} J. H. Keum, 
Automorphisms of Jacobian Kummer surfaces. Compositio Math. 107 (1997), no.
3, 269--288. 
\item{[K-K]} J. H. Keum, S. Kondo, The automorphism  groups of Kummer surfaces
associated with the product of two elliptic curves, 1998 preprint. 
\item{[K]} S. Kondo, The automorphism group of a generic Jacobian Kummer
surface. J. Algebraic Geometry 7 (1998) 589--609.
\item{[K98]} S. Kondo, Niemeier lattices, Mathieu groups and 
finite groups of symplectic automorphisms of K3 surfaces. 
Duke Math. Journal Vol 92 (1998) no. 3, 593-598.
\item{[M]} S. MacLane, Categories for the working mathematician.
Graduate Texts in Mathematics, Vol. 5. Springer-Verlag, New
York-Berlin, 1971. 
\item{[Ma]} G. A. Margulis,
``Discrete subgroups of
semisimple Lie groups''. Ergebnisse der Mathematik und ihrer
Grenzgebiete (3),
17. Springer-Verlag, Berlin, 1991.  ISBN: 3-540-12179-X.
\item{[Mu]} S. Mukai, 
Finite groups of automorphisms of K3 surfaces and the Mathieu group. Invent.
Math. 94 (1988), no. 1, 183--221.
\item{[N]} V. V. Nikulin,  An analogue of the Torelli
theorem for Kummer surfaces of Jacobians. (Russian) Izv. Akad.  Nauk
SSSR Ser. Mat. 38 (1974), 22--41. Integral symmetric bilinear forms
and its applications, Math. USSR Izv., 14 (1980), 103-167.
\item{[P-S]} I. R. Shafarevich, I. I. Piatetski-Shapiro,
Torelli's theorem for algebraic surfaces of type K3.
(Russian) 
Izv. Akad. Nauk SSSR Ser. Mat. 35 1971 530--572.
English translation in pages 516-557 of 
Collected mathematical papers by Igor R. Shafarevich.
Springer-Verlag, Berlin-New York, 1989.  ISBN: 3-540-13618-5 
\item{[Q]} D. Quillen, Higher algebraic $K$-theory. I.
Algebraic $K$-theory, I: Higher $K$-theories (Proc.
Conf., Battelle Memorial Inst., Seattle, Wash., 1972),
pp. 85--147. Lecture Notes in Math., Vol. 341, Springer, Berlin 1973.
\item{[S]} J.-P. Serre,  $\hbox{Cohomologie des groupes discrets}$.
Prospects in mathematics (Proc.
Symp., Princeton Univ., Princeton, N.J., 1970), pp. 77--169.
Ann. of Math. Studies, No. 70, Princeton Univ. Press,
Princeton, N.J., 1971.
\item{[V]} \`E. B. Vinberg,
The two most algebraic $K3$ surfaces. Math. Ann. 265 (1983), no. 1, 1--21.
\item{[V-K]} \`E. B. Vinberg, I. M. Kaplinskaja,
The groups $O\sb{18,1}(Z)$ and $O\sb{19,1}(Z)$. (Russian)
Dokl. Akad. Nauk SSSR 238 (1978), no. 6, 1273--1275.
\bye